\documentclass[final,1p,times]{elsarticle}

\usepackage{amssymb}
 \usepackage{amsthm}
\usepackage{amscd}
\usepackage{amsmath}
\usepackage{amsfonts}
\usepackage{amssymb}
\usepackage{graphicx}
\newtheorem{theorem}{Theorem}[section]

\usepackage{mathrsfs}
\usepackage{titletoc}
\numberwithin{equation}{section}

\newcommand{\D}{\Delta}

\newcommand{\f}{\frac}
\newcommand{\al}{\alpha}

\renewcommand{\l}{\lambda}

\newcommand{\be}{\begin{equation}}

\newcommand{\ee}{\end{equation}}
\newcommand{\bea}{\begin{eqnarray}}
\newcommand{\eea}{\end{eqnarray}}
\newcommand{\bna}{\begin{eqnarray*}}
\newcommand{\ena}{\end{eqnarray*}}
\renewcommand{\o}{\omega}

\renewcommand{\le}{\left}
\newcommand{\ri}{\right}

\newcommand{\na}{\nabla}

\journal{***}
\begin{document}

\begin{frontmatter}

\title{ Gradient estimates for the nonlinear parabolic equation with two exponents on Riemannian manifolds}

\author{Songbo Hou \corref{cor1}}
\ead{housb@cau.edu.cn}
\address{Department of Applied Mathematics, College of Science, China Agricultural University,  Beijing, 100083, P.R. China}

\cortext[cor1]{Corresponding author: Songbo Hou}

\begin{abstract}

In this paper, we study the nonlinear parabolic equation with two exponents on complete noncompact Riemannian maniflods. The special types of such equation include the Fisher-KPP equation, the parabolic Allen-Cahn equation and the Newell-Whitehead equation. We get the Souplet-Zhang's gradeint estimates for the positive solutions to such equation.  We also obtain the Liouville theorem for positive ancient solutions. Our results extend those of Souplet-Zhang (Bull. London. Math. Soc. 38:1045-1053, 2006) and Zhu (Acta Mathematica Scientia 36B(2): 514-526, 2016).
\end{abstract}

\begin{keyword}  Gradient estimate\sep Nonlinear parabolic equation, Liouville theorem.

\MSC [2020] 35K55, 58J35

\end{keyword}

 \end{frontmatter}

\section{Introduction}

Let $M$ be  a complete noncompact Riemannian manifold. In this paper,
 we consider the following nonlinear parabolic equation
\be \label{eq1}\frac{\partial u}{\partial t}=\D u(x,t) +\l(x,t)u^p+\eta (x,t)u^q
\ee on $M$,
 where the functions $\l$ and $\eta$ are $C^1$ in $x$  and $C^0$ in $t$, $p$ and $q$ are positive constants with  $p\geq 1$, $q\geq 1$. If
$\l =-\eta=c$, , $p=1$ and $q=2$, where $c$ is a positive constant, then the equation (\ref{eq1}) becomes
\be \label{kpp}\frac{\partial u}{\partial t}=\D u +cu(1-u)
\ee
which is called the Fisher-KPP equation \cite{Fi37, KPP}.   It describes
the propagation of an evolutionarily advantageous gene in a population and has many applications. Cao,  Liu,   Pendleton and Ward \cite{CLP} derived some  differential Harnack estimates  for positive solutions to (\ref{kpp}) on Riemannian manifolds.  Geng and the author \cite{GH} extended the result of \cite{CLP}. If $\l=1$, $\eta=-1$, $p=1$ and $q=3$, then the equation (\ref{eq1}) becomes
\be
\frac{\partial u}{\partial t}=\D u -(u^3-u)
\ee
which is called the parabolic Allen-Cahn equation. A Harnack inequality for this equation was studied  in \cite{Mhi}. The gradient estimates for the elliptic Allen-Cahn equation on Riemannian manifolds were obtained by the author in \cite {Hou19}.  The special type of (\ref{eq1}) also includes
the Newell-Whitehead equation \cite{NW69}
\be
\frac{\partial u}{\partial t}=\D u +au-bu^3
\ee
where $a$ and $b$ are positive constants. It  is used to model the change of
concentration of a substance. The reader may refer to \cite{BBC} for the recent results for such equation.

The gradient estimate is an important method in study on parabolic and elliptic equations. It was first proved  by Yau \cite{Y75} and Cheng-Yau \cite{CY75}, and was further developed by Li-Yau \cite{LY86},
Li \cite{Li91}, Hamilton \cite{Ha93}, Negrin \cite{N95}, Souplet and Zhang \cite{SZ06}, Ma \cite{Ma06}, Yang \cite{Yang}, etc.
In \cite{SZ06}, Souplelt and Zhang considered the heat equation
\be \label{hte}\f{\partial u}{\partial t}=\D u\ee and proved the following result.

\vskip 6pt

\noindent
{\bf Theorem A.} {\it Let $M$ be  an n-dimensional Riemannian manifold with $n\geq 2$ and $Ricci(M)\geq -k$, $k\geq 0$.
 If  $u$ is any positive solution to (\ref{hte}) in $Q_{R,T}\equiv B(x_0, R)\times [t_0-T, t_0]\subset M\times(-\infty, +\infty)$ and $u\leq N$
 in $Q_{R,T}$, then there holds
 \be\f{|\na u(x,t)|}{u(x,t)}\leq c\le(\f{1}{R}+\f{1}{T^{\f{1}{2}}}+\sqrt{k}\ri)\le(1+\log\f{N}{u(x,t)}\ri)\ee
  in $Q_{\f{R}{2},\f{T}{2}}$, where $c=c(n)$.}

\vskip 6pt
 Later, using the method of Souplet and Zhang,  Zhu \cite{Zhu} studied the equation
 \be\label{npe}
 \le( \D-\f{\partial}{\partial t}\ri)u(x,t)+h(x,t)u^p(x,t)=0,\,\,p>1
 \ee
 on compete noncomapct Riemannian manfolds, where the function $h(x,t)$ is  assumed to be $C^1$ in the first variable and $C^0$ in the second variable. He proved the following result.

\vskip 6pt
  \noindent
{\bf Theorem B.} {\it Let $M$  be an n-dimensional Riemannian manifold with $n\geq 2$ and $Ricci(M)\geq -k$, $k\geq 0$.
 If $u$ is any positive solution to (\ref{npe}) in $Q_{R,T}\equiv B(x_0, R)\times [t_0-T, t_0]\subset M\times(-\infty, +\infty)$ and $u\leq N$
 in $Q_{R,T}$,  then  for any $\beta\in(0,2)$,  there exists a  constant $c=c(n,p,\beta)$ such that
 \be\label{zin}\f{|\na u(x,t)|^2}{u(x,t)^{\beta}}\leq cN^{2-\beta}\le(\f{1}{R^2}+\f{1}{T}+k+N^{p-1}||h^+||_{L^{\infty}(Q_{R,T})}
 +N^{\f{2}{3}(p-1)}||\na h||^{\f{2}{3}}_{L^{\infty}(Q_{R,T})}\ri)\ee
  in $Q_{\f{R}{2},\f{T}{2}}$, where $h^+=\max\{h,0\}$.}

 \vskip 6pt

The same method was also used by Huang and Ma \cite{HM16} to obtain gradient estimates for the equations
$$\f{\partial u}{\partial t}=\D u+\l u^{\al}$$
and $$\f{\partial u}{\partial t}=\D u+au\log u +bu$$ under the Ricci flow, where $\l$,  $\al$, $a$ and $b$ are constants.
\vskip 6pt
  In this paper, we get the following result.
  \begin{theorem}
 Let $M$  be an n-dimensional Riemannian manifold with $n\geq 2$ and $Ricci(M)\geq -k$, $k\geq 0$. Suppose that $\lambda(x,t)$ and $ \eta(x,t)$ are $C^1$ in $x$  and $C^0$ in $t$, $p$ and $q$ are positive constants with  $p\geq 1$, $q\geq 1$.
 If  $u$ is any positive solution to (\ref{eq1}) in $Q_{R,T}\equiv B(x_0, R)\times [t_0-T, t_0]\subset M\times(-\infty, +\infty)$ and $u\leq N$
 in $Q_{R,T}$,  then  there exists a  constant $c=c(n,p,q)$ such that
\be\label{min}\begin{split}
\f{|\na u(x,t)|}{u(x,t)}&\leq c\le( \f{1}{R}+\f{1}{\sqrt{T}}+\sqrt{k}+N^{\f{p-1}{2}}||\l^{+}||^{\f{1}{2}}_{L^{\infty}(Q_{R,T})}
+N^{\f{q-1}{2}}||\eta^{+}||^{\f{1}{2}}_{L^{\infty}(Q_{R,T})}\ri.\\
&\le.\quad+N^{\f{1}{3}(p-1)}||\na \l||^{\f{1}{3}}_{L^{\infty}(Q_{R,T})}+N^{\f{1}{3}(q-1)}||\na \eta||^{\f{1}{3}}_{L^{\infty}(Q_{R,T})}\ri)
\le(1+\log\f{N}{u}\ri)
\end{split}
\ee
  in $Q_{\f{R}{2},\f{T}{2}}$, where $\l^+=\max\{\l,0\}$,  $\eta^+=\max\{\eta,0\}$.
  \end{theorem}

  Note that the estimate (\ref{zin}) is equivalent to
   \be\label{zin2}\f{|\na u(x,t)|}{u(x,t)}\leq c\le(\f{N}{u}\ri)^{1-\f{\beta}{2}}\le(\f{1}{R}+\f{1}{\sqrt{T}}+\sqrt{k}+N^{\f{p-1}{2}}||h^+||^{\f{1}{2}}_{L^{\infty}(Q_{R,T})}
 +N^{\f{1}{3}(p-1)}||\na h||^{\f{1}{3}}_{L^{\infty}(Q_{R,T})}\ri).\ee
 Appling Theorem 1.1 to (\ref{npe}) yields
 \be\label{zin3}\f{|\na u(x,t)|}{u(x,t)}\leq c\le(\f{1}{R}+\f{1}{\sqrt{T}}+\sqrt{k}+N^{\f{p-1}{2}}||h^+||^{\f{1}{2}}_{L^{\infty}(Q_{R,T})}
 +N^{\f{1}{3}(p-1)}||\na h||^{\f{1}{3}}_{L^{\infty}(Q_{R,T})}\ri)\le(1+\log \f{N}{u}\ri).\ee

 Since $\lim\limits_{x\rightarrow +\infty}\f{\log x}{x^{1-\f{\beta}{2}}}=0$, if $\f{N}{u}$ is large enough, then we have
 \be 1+\log\f{N}{u}\leq \le(\f{N}{u}\ri)^{1-\f{\beta}{2}}. \ee  So in this sense, the estimate (\ref{min}) improves (\ref{zin}).

 We also get the Liouville type theorem.

  \begin{theorem}
 Let $M$  be an n-dimensional Riemannian manifold with nonnegative Ricci curvature.
 Suppose that $\l$, $\eta$ are nonpositive constants and one of them is negative,  then equation (\ref{eq1}) does not admit any positive ancient solution with $u(x,t)=e^{o(d(x)+\sqrt{|t|})}$
near infinity.
 \end{theorem}
 The method of the proofs of main theorems comes from \cite{SZ06}, \cite{Zhu} and \cite{HM16}.
 \section{Proof of main theorems}

 \subsection{Proof of Theorem 1.1}

 Let $\tilde{u}=u/N$. Then $\tilde{u}$ satisfies
 \be \label{eq2}\frac{\partial \tilde{u}}{\partial t}=\D \tilde{u} +\tilde{\l}\tilde{u}^p+\tilde{\eta}\tilde{u}^q
\ee
where $\tilde{\l}=\l N^{p-1}$, $\tilde{\eta }=\eta N^{q-1}$. Noting $\tilde{u}\leq 1$, we let
 \be \label{eqo}f=\log \tilde{u}, \,\,\,\,\,\omega=|\na \ln(1-f)|^2.\ee In view of (\ref{eq2}), we have
\be\label {eqf}
\D f +|\na f|^2 +\tilde{\l}e^{(p-1)f}+\tilde{\eta}e^{(q-1)f}-f_t=0.
\ee
By (\ref{eqo}) and (\ref{eqf}), we have
\be \begin{split}
\o_t =&\f{2f_i(f_t)_i}{(1-f)^2}+\f{2f_j^2f_t}{(1-f)^3}\\
= &\f{2f_i\le(f_{jji}+2f_jf_{ji}+\tilde{\l}_ie^{(p-1)f}+\tilde{\l}(p-1)e^{(p-1)f}f_i+\tilde{\eta}_ie^{(q-1)f}+\tilde{\eta}(q-1)e^{(q-1)f}f_i\ri)}{(1-f)^2} \\
&+\f{2f_j^2(f_{ii}+f_i^2+\tilde{\l}e^{(p-1)f}+\tilde{\eta}e^{(q-1)f})}{(1-f)^3}.
\end{split}
\ee
It follows from the similar calculation that
\be
\begin{split}
\D \o&=\f{2f_{ij}^2+2f_jf_{jii}}{(1-f)^2}+\f{8f_if_{ij}f_j+2f_j^2f_{ii}}{(1-f)^3}+\f{6f_i^2f_j^2}{(1-f)^4}\\
&=\f{2f_{ij}^2+2f_jf_{iij}+2R_{ij}f_if_j}{(1-f)^2}+\f{8f_if_{ij}f_j+2f_j^2f_{ii}}{(1-f)^3}+\f{6f_i^2f_j^2}{(1-f)^4}
\end{split}
\ee
where Bochner's identity is used.
 Noting that $R_{ij}f_if_j\geq -kf_i^2$, we have
\be\label{wte}
\begin{split}
\D \o-\omega_t&\geq \f{2f_{ij}^2-4f_if_jf_{ij}-2e^{(p-1)f}f_i\tilde{\l}_i-2\tilde{\l}(p-1)e^{(p-1)f}f_i^2}{(1-f)^2}\\
&\quad-\f{2e^{(q-1)f}f_i\tilde{\eta}_i+2\tilde{\eta}(q-1)e^{(q-1)f}f_i^2+2kf_i^2}{(1-f)^2}\\
&\quad+\f{8f_if_{ij}f_j-2f_j^2f_i^2-2\tilde{\l}e^{(p-1)f}f_j^2-2\tilde{\eta}e^{(q-1)f}f_j^2}{(1-f)^3}+\f{6f_i^2f_j^2}{(1-f)^4}.
\end{split}
\ee
From (\ref{eqo}), we deduce that
\be\label{for}
-\f{2f}{1-f}\na f\na \o=\f{4f_if_{ij}f_j}{(1-f)^2}+\f{4f_i^2f_j^2-4f_if_{ij}f_j}{(1-f)^3}-\f{4f_i^2f_j^2}{(1-f)^4}.
\ee
Combining (\ref{wte}) and (\ref{for}), we have
\be\label{wt2}
\begin{split}
\D \o-\omega_t-\f{2f}{1-f}\na f\na \o&\geq \f{2f_{ij}^2-2e^{(p-1)f}f_i\tilde{\l}_i-2\tilde{\l}(p-1)e^{(p-1)f}f_i^2}{(1-f)^2}\\
&\quad-\f{2e^{(q-1)f}f_i\tilde{\eta}_i+2\tilde{\eta}(q-1)e^{(q-1)f}f_i^2+2kf_i^2}{(1-f)^2}\\
&\quad+\f{4f_if_{ij}f_j+2f_j^2f_i^2-2\tilde{\l}e^{(p-1)f}f_j^2-2\tilde{\eta}e^{(q-1)f}f_j^2}{(1-f)^3}+\f{2f_i^2f_j^2}{(1-f)^4}.
\end{split}
\ee
H\"{o}lder's inequality implies that
\be
\le| \f{4f_if_{ij}f_j}{(1-f)^3}\ri|\leq \f{2f_{ij}^2}{(1-f)^2}+\f{2f_i^2f_j^2}{(1-f)^4}.
\ee
Thus we have
\be\label{wt3}
\begin{split}
\D \o-\omega_t-\f{2f}{1-f}\na f\na \o&\geq -\f{2e^{(p-1)f}f_i\tilde{\l}_i+2e^{(q-1)f}f_i\tilde{\eta}_i}{(1-f)^2}\\
&\quad-\f{2\tilde{\l}(p-1)e^{(p-1)f}f_i^2+2\tilde{\eta}(q-1)e^{(q-1)f}f_i^2+2kf_i^2}{(1-f)^2}\\
&\quad+\f{2f_i^2f_j^2-2\tilde{\l}e^{(p-1)f}f_j^2-2\tilde{\eta}e^{(q-1)f}f_j^2}{(1-f)^3}\\
&=2(1-f)\o^2-2\tilde{\l}\le(p-1+\f{1}{1-f}\ri)e^{(p-1)f}\o\\
&\quad-2\tilde{\eta}\le(q-1+\f{1}{1-f}\ri)e^{(q-1)f}\o\\
&\quad-\f{2e^{(p-1)f}f_i\tilde{\l}_i+2e^{(q-1)f}f_i\tilde{\eta}_i}{(1-f)^2}-2k\o.
\end{split}
\ee
Now we choose a smooth cut-off function $\psi=\psi(x,t)$ with compact support in $Q_{R,T}$ such that

(1) $\psi =\psi (r,t)$, $0\leq \psi\leq 1$ with $\psi=1$ in $Q_{R/2, T/2}$, where $r=d(x,x_0)$;

(2) $\psi$ is decreasing with respect to $r$;

(3) for any $0<\al<1$, $|\partial _r \psi|/\psi^{\al}\leq C_\al/R$, $|\partial_r^2\psi|/\psi^\al\leq C_\al/R^2$ ;

(4) $|\partial_t\psi|/\psi^{1/2}\leq C/T$.

Using (\ref{wt3}), we get
\be\label{wt4}
\begin{split}
\D (\psi\o)-2\frac{\na \psi}{\psi}\cdot \na(\psi\o)-(\psi\omega)_t&\geq 2(1-f)\psi\o^2-2\tilde{\l}\le(p-1+\f{1}{1-f}\ri)e^{(p-1)f}\psi\o\\
&\quad-2\tilde{\eta}\le(q-1+\f{1}{1-f}\ri)e^{(q-1)f}\psi\o\\
&\quad-2k\psi\o-\f{2e^{(p-1)f}f_i\tilde{\l}_i+2e^{(q-1)f}f_i\tilde{\eta}_i}{(1-f)^2}\psi\\
&\quad+\f{2f}{1-f}\na f\na(\psi \o)-\f{2f\o}{1-f}\na f\na\psi-\f{2|\na \psi|^2}{\psi}\o\\
&\quad+(\D \psi)\o-\psi_t\o.
\end{split}
\ee
 Suppose that $\psi\o$ attains the maximum at $(x_1, t_1)$.  The argument in \cite{Ca58} implies that  we can assume  $x_1$ is not in the cut-locus of $M$. Then we have $\D(\psi \o)\leq 0$, $(\psi\o)_t\geq 0$ and $\na(\psi \o)=0$ at $(x_1,t_1)$.
 It follows that
 \be\label{wt5}
\begin{split}
 2(1-f)\psi\o^2 &\leq    2\tilde{\l}\le(p-1+\f{1}{1-f}\ri)e^{(p-1)f}\psi\o +2\tilde{\eta}\le(q-1+\f{1}{1-f}\ri)e^{(q-1)f}\psi\o\\
&\quad+2k\psi\o+\f{2e^{(p-1)f}f_i\tilde{\l}_i+2e^{(q-1)f}f_i\tilde{\eta}_i}{(1-f)^2}\psi\\
&\quad+\f{2f\o}{1-f}\na f\na\psi+\f{2|\na \psi|^2}{\psi}\o-(\D \psi)\o+\psi_t\o.
\end{split}
\ee
In view of $p\geq 1$, $q\geq 1$ and $f\leq 0$,   we have
 \be\label{dpe}
\begin{split}
&2\tilde{\l}\le(p-1+\f{1}{1-f}\ri)e^{(p-1)f}\psi\o +2\tilde{\eta}\le(q-1+\f{1}{1-f}\ri)e^{(q-1)f}\psi\o\\
&\quad\leq 2\tilde{\l}^{+}p\psi\o+2\tilde{\eta}^{+}q\psi\o\\
&\quad\leq \f{1}{16}\psi\o^2+16\psi(\tilde{\l}^{+}p)^2+\f{1}{16}\psi\o^2+16\psi(\tilde{\eta}^{+}q)^2\\
&\quad\leq \f{1}{8}\psi\o^2+16(\tilde{\l}^{+}p)^2+16(\tilde{\eta}^{+}q)^2
\end{split}
\ee
where $\tilde{\l}^{+}=\max\{\tilde{\l},0\}$, $\tilde{\eta}^{+}=\max\{\tilde{\eta},0\}$. Straightforward calculations show

\be\label{fie}
\begin{split}
\f{2e^{(p-1)f}f_i\tilde{\l}_i+2e^{(q-1)f}f_i\tilde{\eta}_i}{(1-f)^2}\psi &\leq \f{f_i^4}{2(1-f)^4}\psi+\f{3|\na\tilde{\l}|^{4/3}}{2(1-f)^{4/3}}\psi+\f{f_i^4}{2(1-f)^4}\psi+\f{3|\na\tilde{\eta}|^{4/3}}{2(1-f)^{4/3}}\psi\\
&\leq \f{f_i^4}{(1-f)^4}\psi+\f{3}{2}(|\na\tilde{\l}|^{4/3}+|\na\tilde{\eta}|^{4/3})\\
&\leq  (1-f)\psi\o^2+\f{3}{2}(|\na\tilde{\l}|^{4/3}+|\na\tilde{\eta}|^{4/3}),
\end{split}
\ee
\be\label{bfe}
\begin{split}
\le|\f{2f\o}{1-f}\na f\na\psi\ri|&\leq 2\o^{3/2}|f||\na \psi|=2[\psi(1-f)\o^2]^{3/4}\f{|f||\na \psi|}{[\psi(1-f)]^{3/4}}\\
&\leq \f{1}{8} (1-f)\psi\o^2+c\f{(f|\na\psi|)^4}{[\psi(1-f)]^3}\\
&\leq \f{1}{8} (1-f)\psi\o^2+c\f{f^4}{R^4(1-f)^3},
\end{split}
\ee
\be\label{koe}
\begin{split}
2k\psi\o\leq \f{1}{8}(1-f)\psi\o^2+ck^2.
\end{split}
\ee
By the estimates of Souplet and Zhang \cite{SZ06}, we have
\be\label{pse}
\f{|\na \psi|^2}{\psi}\o\leq \f{1}{8}\psi\o^2+c\f{1}{R^4}\leq \f{1}{8}(1-f)\psi\o^2+c\f{1}{R^4},
\ee
\be\label{dpe}
-(\D\psi)\o\leq \f{1}{8}\psi\o^2+c\f{1}{R^4}+ck\f{1}{R^2}\leq \f{1}{8}(1-f)\psi\o^2+c\f{1}{R^4}+ck\f{1}{R^2},
\ee
\be\label{pte}
|\psi_t|\o\leq \f{1}{8}\psi \o^2+c\f{1}{T^2}\leq \f{1}{8}(1-f)\psi \o^2+c\f{1}{T^2} .
\ee
Combining (\ref{wt5})-(\ref{pte}), we obtain
\be
\begin{split}
\f{1}{8}(1-f)\psi\o^2&\leq 16(\tilde{\l}^{+}p)^2+16(\tilde{\eta}^{+}q)^2+\f{3}{2}(|\na\tilde{\l}|^{4/3}+|\na\tilde{\eta}|^{4/3})\\
&\quad+c\f{f^4}{R^4(1-f)^3}+ck^2+c\f{1}{R^4}+ck\f{1}{R^2}+c\f{1}{T^2}.
\end{split}
\ee

Hence
\be
\begin{split}
\psi \o^2(x_1,t_1)&\leq  cN^{2p-2}||\l^{+}||^2_{L^{\infty}(Q_{R,T})}+cN^{2q-2}||\eta^{+}||^2_{L^{\infty}(Q_{R,T})}\\
&\quad+cN^{\f{4}{3}(p-1)}||\na \l||^{\f{4}{3}}_{L^{\infty}(Q_{R,T})}+cN^{\f{4}{3}(q-1)}||\na \eta||^{\f{4}{3}}_{L^{\infty}(Q_{R,T})}\\
&\quad+c\f{f^4}{R^4(1-f)^4}+ck^2+c\f{1}{R^4}+c\f{1}{T^2}.
\end{split}
\ee
By above estimate, there holds for all $(x,t)$ in $Q_{R,T}$,
\be
\begin{split}
\psi^2 \o^2(x,t)&\leq  cN^{2p-2}||\l^{+}||^2_{L^{\infty}(Q_{R,T})}+cN^{2q-2}||\eta^{+}||^2_{L^{\infty}(Q_{R,T})}\\
&\quad+cN^{\f{4}{3}(p-1)}||\na \l||^{\f{4}{3}}_{L^{\infty}(Q_{R,T})}+cN^{\f{4}{3}(q-1)}||\na \eta||^{\f{4}{3}}_{L^{\infty}(Q_{R,T})}\\
&\quad+c\f{1}{R^4}+c\f{1}{T^2}+ck^2.
\end{split}
\ee
Noting that $\psi(x,t)=1$ in $Q_{R/2, T/2}$, we get
\be\begin{split}
\f{|\na f(x,t)|}{1-f(x,t)}&\leq \f{c}{R}+\f{c}{\sqrt{T}}+c\sqrt{k}+cN^{\f{p-1}{2}}||\l^{+}||^{\f{1}{2}}_{L^{\infty}(Q_{R,T})}
+cN^{\f{q-1}{2}}||\eta^{+}||^{\f{1}{2}}_{L^{\infty}(Q_{R,T})}\\
&\quad+cN^{\f{1}{3}(p-1)}||\na \l||^{\f{1}{3}}_{L^{\infty}(Q_{R,T})}+cN^{\f{1}{3}(q-1)}||\na \eta||^{\f{1}{3}}_{L^{\infty}(Q_{R,T})}.
\end{split}
\ee
Finally we have
\be\begin{split}
\f{|\na u(x,t)|}{u(x,t)}&\leq  c\le( \f{1}{R}+\f{1}{\sqrt{T}}+\sqrt{k}+N^{\f{p-1}{2}}||\l^{+}||^{\f{1}{2}}_{L^{\infty}(Q_{R,T})}
+N^{\f{q-1}{2}}||\eta^{+}||^{\f{1}{2}}_{L^{\infty}(Q_{R,T})}\ri.\\
&\quad\le.+N^{\f{1}{3}(p-1)}||\na \l||^{\f{1}{3}}_{L^{\infty}(Q_{R,T})}+N^{\f{1}{3}(q-1)}||\na \eta||^{\f{1}{3}}_{L^{\infty}(Q_{R,T})}\ri)
\le(1+\log\f{N}{u}\ri).
\end{split}
\ee
\subsection {Proof of Theorem 1.2}
We prove it by contradiction. Suppose that $u$ is a positive solution to (\ref{eq1}).
Noting that  $\l$ and $\eta$ are nonpositive constants, if follows from Theorem 1.1 that

\be\label{inl}
\f{|\na u(x,t)|}{u(x,t)}\leq c\le( \f{1}{R}+\f{1}{\sqrt{T}}\ri)
\le(1+\log\f{N}{u}\ri).
\ee
By the same argument as in the proof of Theorem 1.2 in \cite{SZ06} and Theorem 1.8 in \cite{Zhu},  fixing $(x_0,t_0)$ and applying (\ref{inl}) to $u$ on $B(x_0,R)\times[t_0-R^2,t_0]$, we get
$$\f{|\na u(x_0,t_0)|}{u(x_0, t_0)}\leq \f{C}{R}[1+o(R)].$$

It follows that $|\na u(x_0,t_0)|=0$ by letting $R\rightarrow\infty$.  Noting $(x_0,t_0)$ is arbitrary, we have $u(x,t)=u(t)$. Then by (\ref{eq1}), we get
$\f{du}{dt}=\l u^p+\eta u^q$. Without loss of generality, we assume that $\l<0$.

If $p>1$, integrating $\f{du}{dt}$ on $[t,0]$ with $t<0$ implies that
$$\f{1}{1-p}(u^{1-p}(0)-u^{1-p}(t))\leq -\l t.$$ Then
$$u^{p-1}(t)\leq u^{p-1}(0)+(1-p)\l t.$$  This yields that if $t$ is large enough, $u^{p-1}(t)<0$ which contradicts that $u$  is positive.

If $p=1$, we get for $t<0$
$$\log u(0)-\log u(t)\leq -\l t.$$
Hence $u(t)\geq u(0)e^{\l t}$, which contradicts $u(x,t)=e^{o(d(x)+\sqrt{|t|})}$
near infinity. We finish the proof.

\end{document}